\documentclass[12pt,a4paper]{article}
\usepackage{kerkis} 
\usepackage{graphics}
\usepackage{amsmath}
\usepackage{amstext}
\usepackage{amssymb}
\usepackage{physics}
\usepackage{color}
\usepackage[scale=0.8]{geometry}

\usepackage[svgnames,dvipsnames]{xcolor}

\begin{document}

\begin{center}
	
	{{{\Large\bf Energy Transport in 1-Dimensional Oscillator Arrays With Hysteretic Damping}}}\\
	
	\vspace{0.5cm}
	{Tassos Bountis\footnote{E-mail address: tassosbountis@gmail.com}$^{,a}$, Konstantinos Kaloudis\footnote{E-mail address: konst.kaloudis@gmail.com}$^{,b}$, Joniald Shena\footnote{E-mail address: jonialdshena@gmail.com}$^{,c}$, Charalampos Skokos\footnote{Corresponding author, e-mail address: haris.skokos@gmail.com}$^{,d}$ and Christos Spitas\footnote{E-mail address: christos.spitas@nu.edu.kz}$^{,b}$}
	
	\vspace{0.2cm}

\footnotesize \phantom{0}\textsuperscript{a} Center for Integrable Systems, P.G. Demidov Yaroslavl State University, Yaroslavl, Russia \newline
\phantom{0}\textsuperscript{b} Department of Mechanical and Aerospace Engineering, Nazarbayev University, Astana, Kazakhstan \newline
\phantom{0}\textsuperscript{c} National University of Science and Technology ``MISiS'', Moscow, Russia\newline
\phantom{0}\textsuperscript{d} Nonlinear Dynamics and Chaos Group, Department of Mathematics and Applied Mathematics, University of Cape Town, Cape Town, South Africa	
\vspace{0.2in}
\end{center}

\begin{abstract}
	Energy transport in 1-dimensional oscillator arrays has been extensively studied to date in the conservative case, as well as under weak viscous damping. When driven at one end by a sinusoidal force, such arrays are known to exhibit the phenomenon of {\it supratransmission}, i.e.~a sudden energy surge above a critical driving amplitude. In this paper, we study 1-dimensional oscillator chains  in the presence of {\it hysteretic damping}, and include nonlinear stiffness forces that are important for many materials at high energies. We first employ Reid's model of {\it local} hysteretic damping, and then study a new model of {\it nearest neighbor dependent} hysteretic damping to compare their supratransmission and wave packet spreading properties in a deterministic as well as stochastic setting. The results have important quantitative differences, which should be helpful when comparing the merits of the two models in specific engineering applications.\\
	
\end{abstract}

\section{Introduction} \label{intro}
In most studies of energy transport in engineering materials, not only particle interactions but also damping and noise phenomena are addressed, while nonlinear characteristics of the materials are also seriously taken into account \cite{OlejAwre2018}. For example, power law terms are introduced to describe more accurately stiffness behavior in the continuum \cite{scalerandi2016power}, while in the case of springs and composite meta-materials quadratic and cubic terms are often proposed \cite{tang2016using,hu2017nonlinear}. In this context, damping mechanisms become crucial and must be carefully modeled to agree with experimental observations \cite{dall,xu}. For example, the familiar choice of {\it viscous} damping \cite{adhikari1} (linearly dependent on the velocity) is often bypassed, in favor of velocity dependent {\it hysteretic} damping, which is deemed more appropriate for oscillating materials whose energy loss per cycle is known to be independent of the deformation frequency \cite{adhikari2,inaudi1995}.

The ultimate purpose of this paper is to contribute new models and results that can be compared with a plethora of theoretical and experimental findings available in the recent literature: For example, in \cite{CS} hysteretic non-causal models of vibration and wave analysis are discussed, \cite{N} proposes a simple hysteretic damping setup is proposed that satisfies a causality condition, while in \cite{CDS} a novel elastic hysteretic model is proposed and analyzed in detail. On the other hand, in \cite{MCRD} hysteretic damping is used to explain experimental measurements of the complex Young’s modulus, while \cite{ZHP} investigates the transverse dynamic hysteretic damping characteristics of a serpentine belt.

In \cite{bks2020}, a detailed study was performed on the dynamics of a hysteretic damping model introduced by Reid \cite{reid1956free}, for a linear as well as nonlinear oscillator driven by a $T$-periodic sinusoidal force. Reid's model was shown to be free from numerical limitations suffered by a hysteretic damping model proposed earlier by Bishop \cite{bishop1955treatment}, since Reid's oscillator is expressed by a differential equation whose real and imaginary parts have the same solutions. Thus, numerical errors occurring during the integration of Reid's model are efficiently controlled, and it is shown to possess $T$-periodic solutions, which are true attractors. Moreover, in the weakly dissipative case, it exhibits the coexistence of periodic attractors of period $nT,n=2,3,...$, with very interesting basins of attraction in the space of initial conditions \cite{bks2020}.

Reid's model \cite{reid1956free} describes the evolution of a one degree--of--freedom (dof) oscillator with mass $M$ under periodic forcing of amplitude $f$, satisfying the differential equation:
\begin{equation}\label{reid1}
M\ddot{x} + c \abs{\dfrac{x}{\dot{x}}} \dot{x}  + k x = M\ddot{x} + k x \left(1 + \dfrac{c}{k} \text{sgn}(x \dot{x})\right) =  f \sin{\omega t},
\end{equation}
where $\text{sgn}(\cdot)$ is the sign function, $x$ denotes the particle's displacement from equilibrium, $c$ is the damping parameter, and $k$ the linear stiffness coefficient.

Eq. (\ref{reid1}) is frequency independent, with work per cycle proportional to the squared amplitude of the oscillations and possesses stable periodic orbits with frequency equal to $\omega=2\pi/T$. Let us now examine the effect of including in our model a nonlinear stiffness term in the form of a cubic term (associated e.g. with a symmetric quartic potential), in which case the equation takes the form:
\begin{equation}\label{reid2}
M\ddot{x} + c\,x\,\text{sgn}(x \dot{x})  + k\,x + \epsilon\,x^3 = M\ddot{x} + c \,\abs{\dfrac{x}{\dot{x}}} \dot{x}  + k \,x +  \epsilon \,x^3 =  f \sin{\omega t}.
\end{equation}

Although the term $ \abs{x/\dot{x}} $ in the original model of Eq. (\ref{reid1}) is nonlinear, this type of nonlinearity does not seriously affect the dynamical behavior of the system, which is often characterized as ``quasilinear'' \cite{elliott}.

In what follows, we avoid any undesired behavior associated with the computation of the sign function in Eq. (\ref{reid2}) using the approximation $\mathrm{sgn}(x-x_0)\approx\tanh[\tau (x-x_0)]$, which works very well for large values of $\tau$ (e.g. $\tau=1000$), as the two functions match for $\tau\to\infty$. We, therefore, write the equation of motion for a single oscillator as
\begin{equation}\label{reid3}
M\ddot{x} + c \abs{x} \tanh (\tau \dot{x}) + k x + \epsilon x^3 =  f \sin{\omega t}.
\end{equation}

We thus begin by studying collective dynamical phenomena arising in 1-dimensional arrays of systems of $N$ nonlinear oscillators of equal mass $M$ under nearest-neighbor coupling, and solve the equations of motion in the presence of Reid's hysteretic damping, as in Eq. (\ref{reid2}). Using the sinusoidal term in the equation as a driving force applied to the first oscillator, we perform numerical computations varying the driving amplitude $f$ and the frequency $\omega$ to obtain a first picture of the collective dynamical behavior of the system under study. 

Our primary objective in this paper is to investigate the occurrence of the important phenomenon of \textit{supratransmission} \cite{geniet2002energy} in two different models with hysteretic damping, as the driving amplitude is gradually increased. As is well--known, supratransmission has been observed in many papers on analytic Hamiltonian models of 1-dimensional oscillators, as a sudden surge of energy in the form of nonlinear waves propagating through the lattice at specific critical amplitudes, over frequency ranges belonging to the forbidden band (see e.g. \cite{maciasbountis2018,bountis2019energy,macias2021nonlinear}). 

Model I is a nonlinear extension of the well-known quasi--linear Reid's model \cite{reid1956free,bks2020}. In Section 2, we numerically integrate its equations of motion and study the occurrence of supratransmission, analyzing the effect of the system's parameters on the corresponding critical amplitudes. We then compare in Section 3 the supratransmission and wave packet spreading properties of model I, with those of a new model, defined here as model II, in which hysteretic damping at each oscillator is {\it nonlocal}, as it involves the positions and velocities of its nearest neighbors. Next, in Section 4, we explore the probability of the occurrence of supratransmission in model I, under {\it random} periodic forcing, where the driving is subjected to stochastic perturbations. 

Next, we have introduced a discussion concerning the phenomenon of {\it breather arrest} that has recently come to our attention \cite{Breather1,Breather2,Breather3,Breather4}. It concerns the propagation of oscillating localized excitations generated impulsively at one end, in lattices of oscillators with damping effects. In Section 5, we describe results indicating that this phenomenon appears, in a similar way, in both our hysteretic damping Models I and II, but with some important differences favoring our hysteretic non-local Model II. Finally, we conclude in Section 6 with a summary and discussion of our results.

\section{Supratransmission in Reid's Hysteretic Damping Model}

In order to derive the equations of motion for the array of $N$ coupled nonlinear Reid oscillators, we begin with a nearest neighbor potential involving quadratic and quartic terms of the form
\begin{equation}\label{pot}
\mathcal{V}\left(\mathbf{x}\right) = \sum_{j=0}^{N} \left[ \dfrac{k}{2}\left(x_{j+1}-x_j\right)^2 + \dfrac{\epsilon}{4}\left(x_{j+1}-x_j\right)^4 \right], 
\end{equation}
which we shall hereafter call model I. Then, using Eqs.~(\ref{reid3}) and (\ref{pot}) with the approximation $\mathrm{sgn}(x)\approx\tanh (\tau x)$ we derive the following equations of motion for $j=1,\ldots,N$ for this model:
\begin{equation}\label{eqm}
M\ddot{x}_j = - c \abs{x_j} \tanh ( \tau \dot{x}_j ) - k\, \left(-x_{j-1} + 2x_j - x_{j+1}\right) - \epsilon \, \left(- \left(x_{j+1}-x_j\right)^3 +\left(x_j-x_{j-1}\right)^3 \right),
\end{equation}
subject to the boundary conditions:
\begin{equation}\label{eq:bc}
x_0(t) = f\sin{\omega\,t} \,\,\,\,\, \text{ and } \,\,\,\,\, x_{N+1}(t)=0,\,\,\,\, t\in \mathcal{T}\subseteq\mathbb{R}^+.
\end{equation}

The mass of all oscillators is set $M=1$ throughout this paper, while for all simulations presented here, we have used the Tsitouras 5/4 Runge-Kutta method \cite{tsit5} with adaptive time step to perform the integration, employing the DifferentialEquations.jl library \cite{rackauckas2017differentialequations} of the Julia programming language \cite{bezanson2017julia}. 
 
As a first step, let us investigate the occurrence of supratransmission in arrays of coupled nonlinear Reid oscillators. To this end, we first need to avoid the so--called forbidden bandgap of the linear mode frequencies \cite{geniet2002energy,maciasbountis2018}, so that they are not excited by our forcing. Since our lattices possess quadratic nearest neighbor terms for $\epsilon=0$ [see Eq. (\ref{pot})], the linear modes are found within $0<\omega<2$. We thus choose from now on to excite our lattices with periodic driving whose frequency satisfies $\omega>2$. Furthermore, since our systems are (hysteretically) damped, traveling waves barely reach the right end of the lattice, which thus acts as an absorbing boundary.

We now apply forcing amplitudes within the interval $f \in \left\{0.2,\ldots, 3.0\right\}$ and keep all other parameters fixed at  $\left(c,k,\omega,\epsilon \right) = \left(0.01, 0.3, 2.5, 0.1\right)$. Starting with all $N=200$ particles at rest, we increase the value of $f$ and plot the time evolution of the array displacements looking for critical amplitudes at which supratransmission occurs. 

\begin{figure}[h!]
	\centering
			\resizebox{0.9\columnwidth}{!}{
		\includegraphics{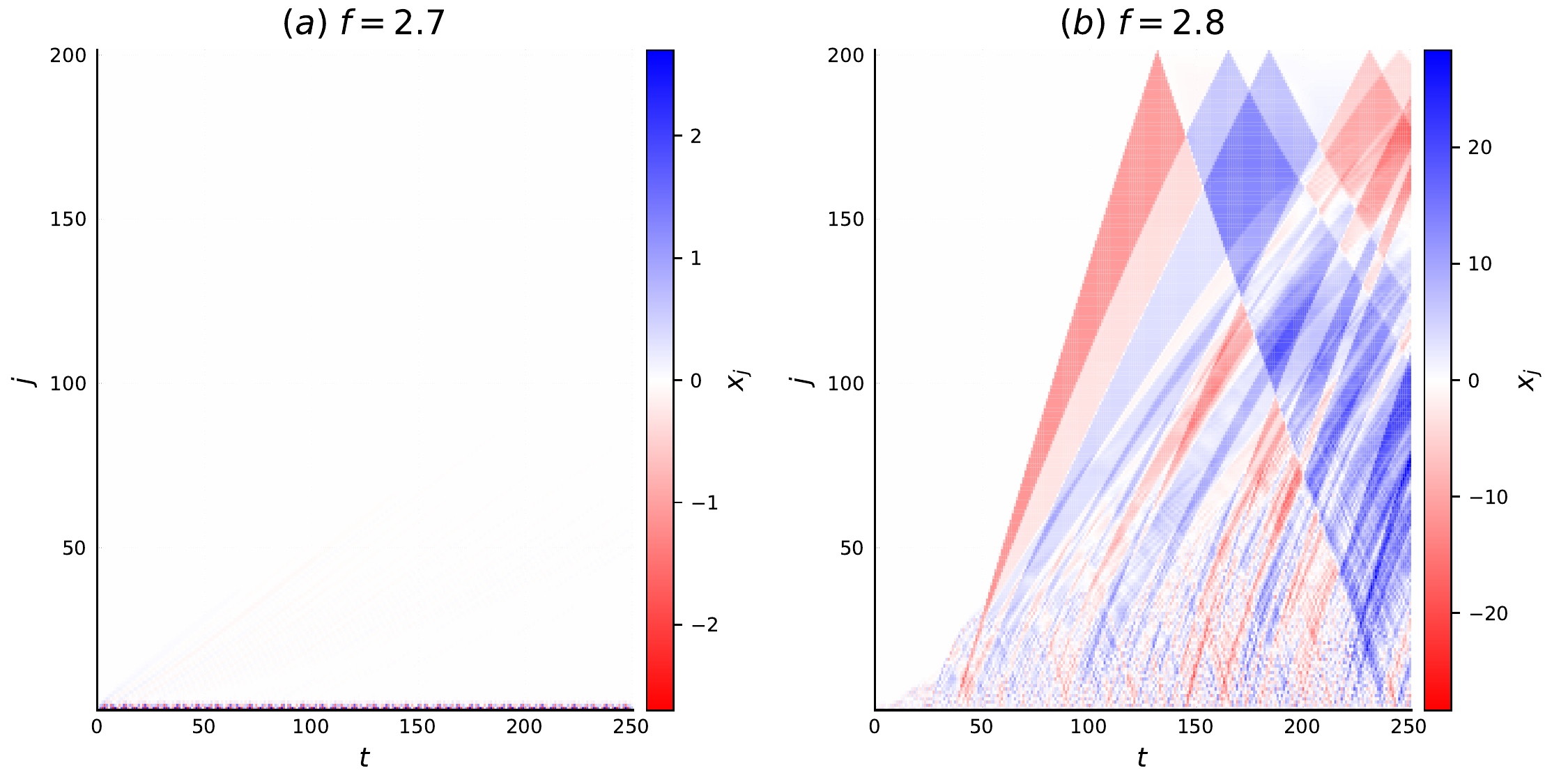} }
	\caption{Graphs of the displacement solutions of Eqs.~(\ref{eqm}) for parameter values $\left(c,k,\omega,\epsilon \right) = \left(0.01, 0.3, 2.5, 0.1\right)$ at (a) $f=2.7$ and (b) $f=2.8$ imply that a large energy wave is generated between these two amplitudes. \label{fig5}}
\end{figure}

Thus, we observe in Fig. \ref{fig5} that, while at $f=2.7$ there is very little energy transfer through the lattice, at $f=2.8$ there is a sudden surge of high amplitude waves, signaling the onset of supratransmission at a critical amplitude lying between these two values. We conjecture that this is due to the excitement of an unstable nonlinear wave mode of the system, as pointed out in earlier studies of Hamiltonian models, where such modes were derived and studied analytically \cite{geniet2002energy,dauxois2007modulational}.

\subsection{Variation of parameters}

In order to better study the onset of supratransmission in the $(\omega, f)$--parameter space of our system, we developed a numerical criterion which identifies more precisely the occurrence of this phenomenon. Specifically, we monitor over time the squared sum of particle displacements and average them, calculating the quantity:
\begin{equation}\label{average} 
\mathcal{D}_f = \dfrac{1}{n_T}\sum_{i=1}^{N}\sum_{j=1}^{n_T}x_i^2(t_j),\,\,\, 0\leq t_1 \leq t_2 ...\leq t_{n_T}
\end{equation}
for various values of the forcing amplitude $f$, keeping all other parameters fixed. We thus identify the critical amplitude $f_{cr}$ as the \textit{first} amplitude whose difference from the previous one exceeds a certain threshold $\Delta_{th}$ (e.g. $\Delta_{th} = 1000$). That is, for a sequence of increasing driving amplitudes $f_1,f_2,\ldots,f_k$ we define as:
\begin{equation}\label{critical}
f_{cr} =  \arg\min_{1\leq i \leq k}\left\{\mathcal{D}_{f_i} - \mathcal{D}_{f_{i-1}} > \Delta_{th}\right\}.
\end{equation}
where the $\mathcal{D}_{f_{i}}$ are calculated using (\ref{average}). Notice that $x_i(\cdot),\, i=1,\ldots,N$ are the solutions of the system of Ordinary Differential Equations (ODEs) (\ref{eqm}), evaluated at discrete time points obtained by numerical integration.

\begin{figure}[h!]
	\centering
			\resizebox{0.8\columnwidth}{!}{
	\includegraphics{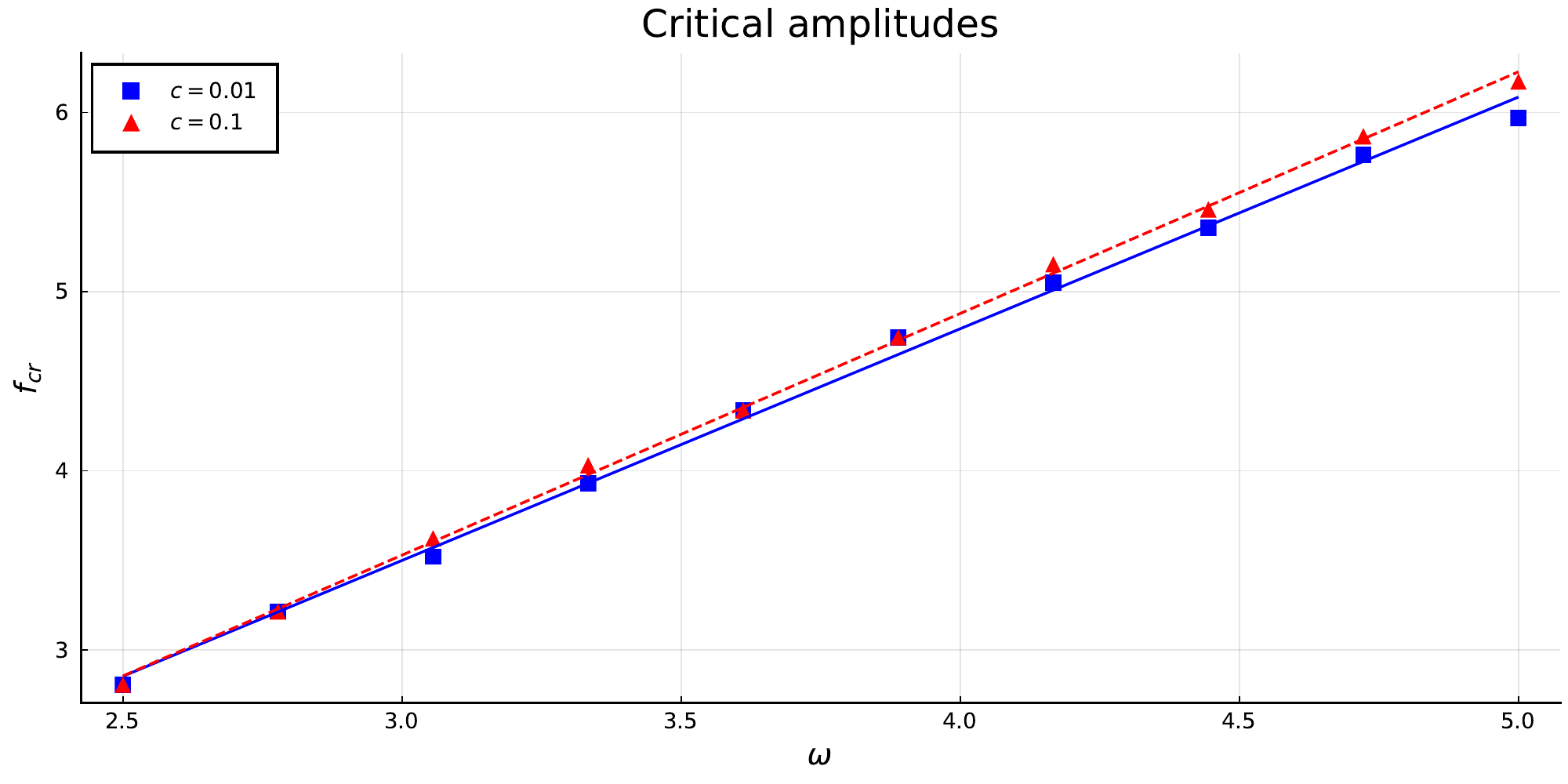} }	
	\caption{Critical forcing amplitudes $f_{cr}$ of model I as functions of the forcing frequency $\omega$, with $\left(k,\epsilon \right) = \left(0.3, 0.1\right)$, for damping parameter $c \in \left\{0.01,0.1\right\}$, and lattices of $N=200$ oscillators. We superimpose the (least squares) linear fits, with slopes 1.29 and 1.34 for $c=0.01$ and $c=0.1$, respectively.}  \label{fig7}
\end{figure}

Let us now study the effect of hysteretic damping on the critical amplitudes. To this end, we fix $\left(k,\epsilon \right) = \left(0.3, 0.1\right)$ and calculate $f_{cr}$ as a function of $\omega$ using the above numerical criterion [see Eqs.~(\ref{average}), (\ref{critical})]. The results presented in Fig. \ref{fig7} evidently show that the damping parameter values do not significantly affect the thresholds $f_{cr}$. In fact, the critical amplitudes increase approximately linearly as functions of $\omega \in \left(2.5,5\right)$, with slopes that vary very little for different damping parameters. 

\begin{figure}[h!]
	\centering
			\resizebox{0.8\columnwidth}{!}{
	\includegraphics{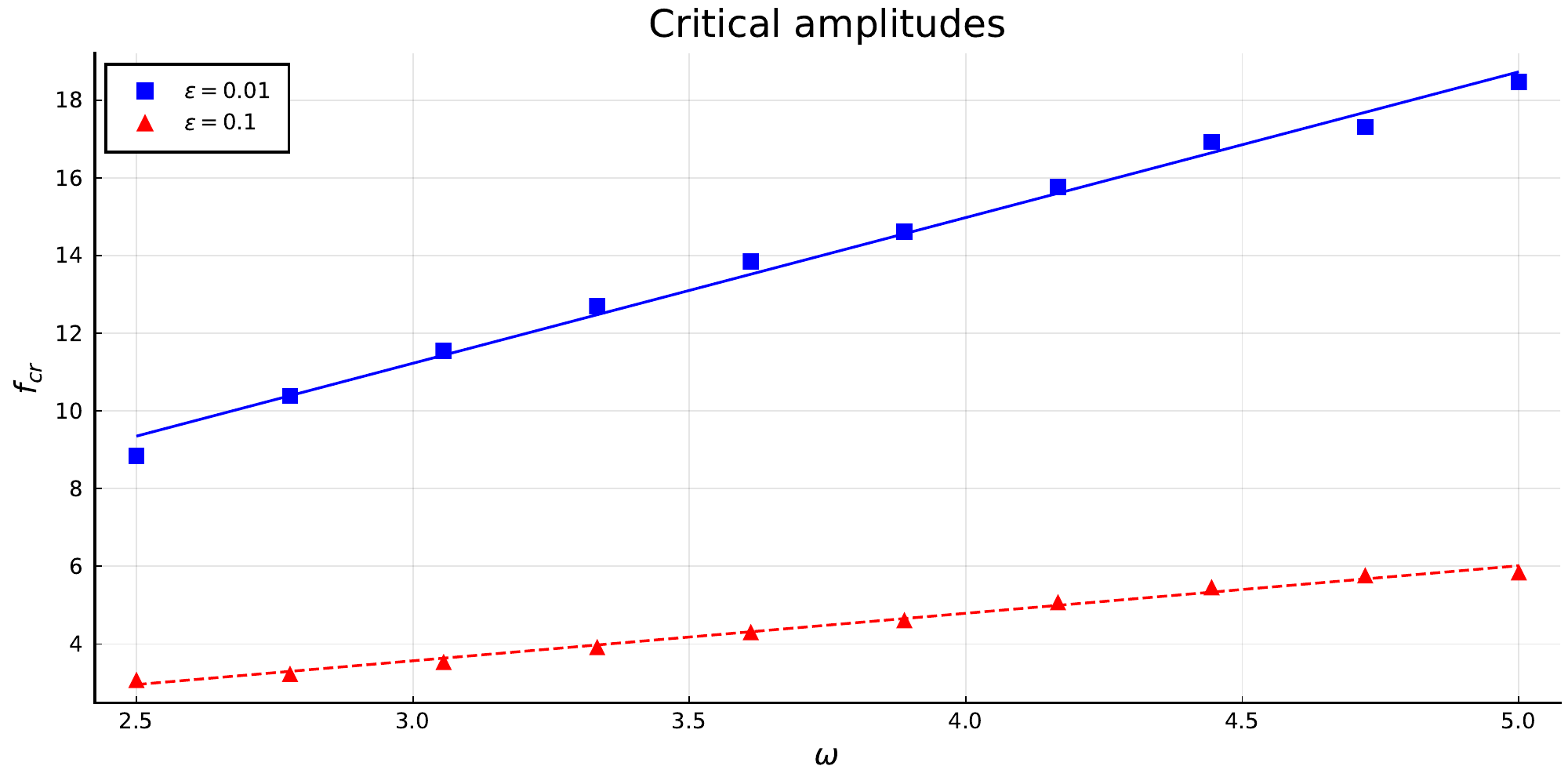} }		
	\caption{Critical forcing amplitudes $f_{cr}$ of model I as functions of the forcing frequency $\omega$, with $\left(k,c \right) = \left(0.3, 0.01\right)$, for nonlinearity parameter $\epsilon \in \left\{0.01,0.1\right\}$, and lattices of $N=200$ oscillators. The dependence on $\omega$ is again linear, with slopes 3.75 and 1.22 for $\epsilon=0.01$ and $\epsilon=0.1$, respectively.\label{fig9}}
\end{figure}

Next, let us examine the effect of varying the nonlinearity parameter $\epsilon$ on the critical amplitudes. To this end, we fix $\left(k,c \right) = \left(0.3, 0.01\right)$ and calculate $f_{cr}$ using the same numerical approach as above. As we observe in Fig. \ref{fig9}, while the dependence of $f_{cr}$ on $\omega$ is still linear, the nonlinearity parameter does indeed have a significant effect, as higher values of $\epsilon$ lead to lower values of $f_{cr}$ as functions of $\omega$. Essentially, this means that when the nonlinear coupling is stronger, it is ``easier'' to excite the nonlinear modes required to give rise to the phenomenon of supratransmission. 

\begin{figure}[h!]
	\centering
			\resizebox{0.8\columnwidth}{!}{
	\includegraphics{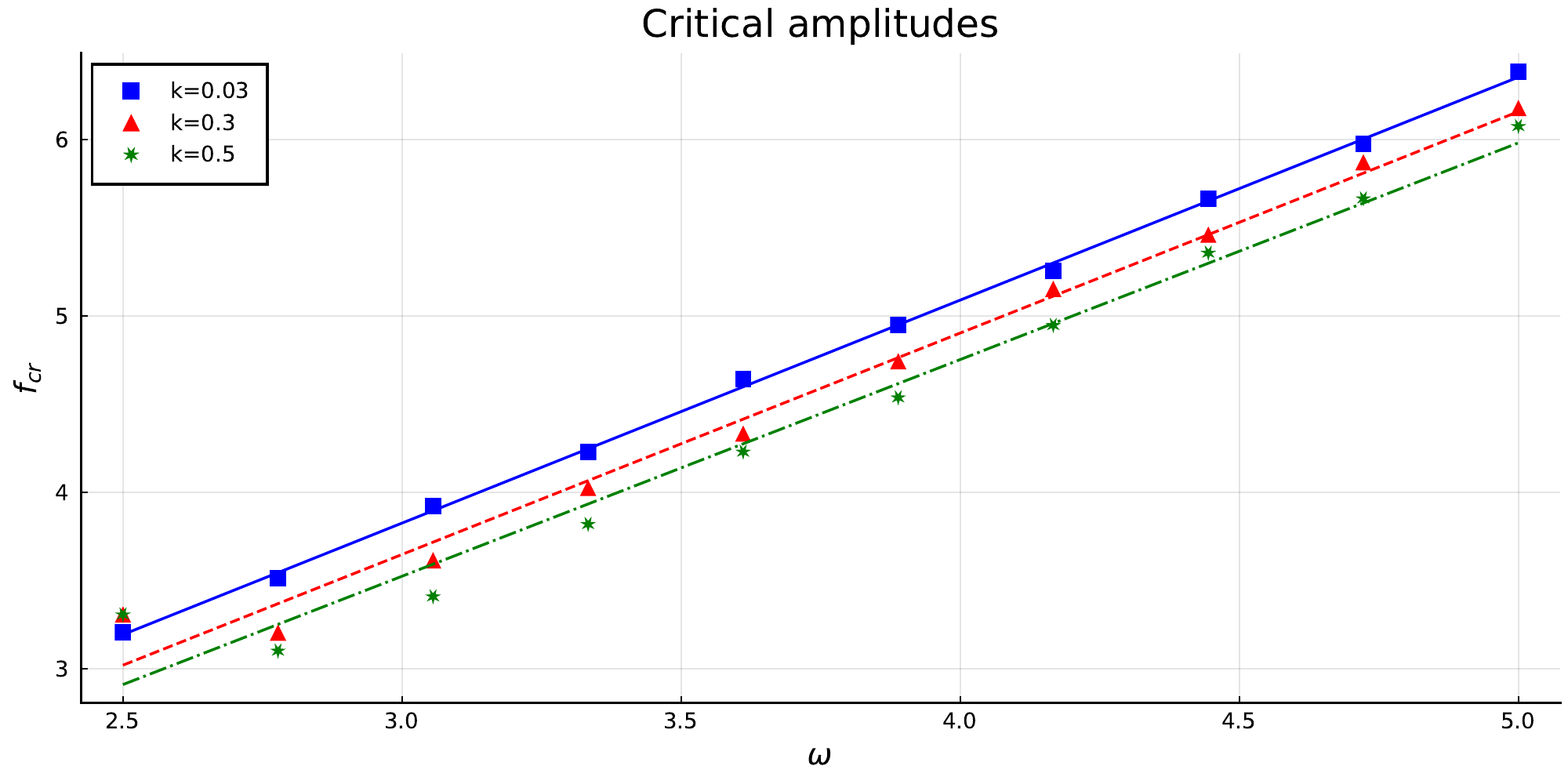} }			
	\caption{Critical forcing amplitudes $f_{cr}$ of model I with respect to the forcing frequency $\omega$, with $\left(c,\epsilon \right) = \left(0.1, 0.1\right)$, for coupling strength $k \in \left\{0.03,0.3, 0.5\right\}$, and lattices of $N=200$ oscillators. We superimpose the (least squares) linear fits, with slopes 1.26, 1.25 and 1.22 for $k=0.03$, $k=0.3$ and $k=0.5$, respectively.\label{fig11}}
\end{figure}

Let us now examine the effect of the linear coupling on the critical amplitudes. To this end, we fix $\left(c,\epsilon \right) = \left(0.1, 0.1\right)$ and calculate $f_{cr}$ as a function of $\omega$ for 3 different $k$ values. The results presented in Fig. \ref{fig11} show that the coupling strength has a similar (but much weaker) effect on $f_{cr}$ as $\epsilon$, while its dependence on $\omega$ is linear with practically the same slope.

\begin{figure}[h!]
	\centering
	\resizebox{0.65\columnwidth}{!}{
		\includegraphics{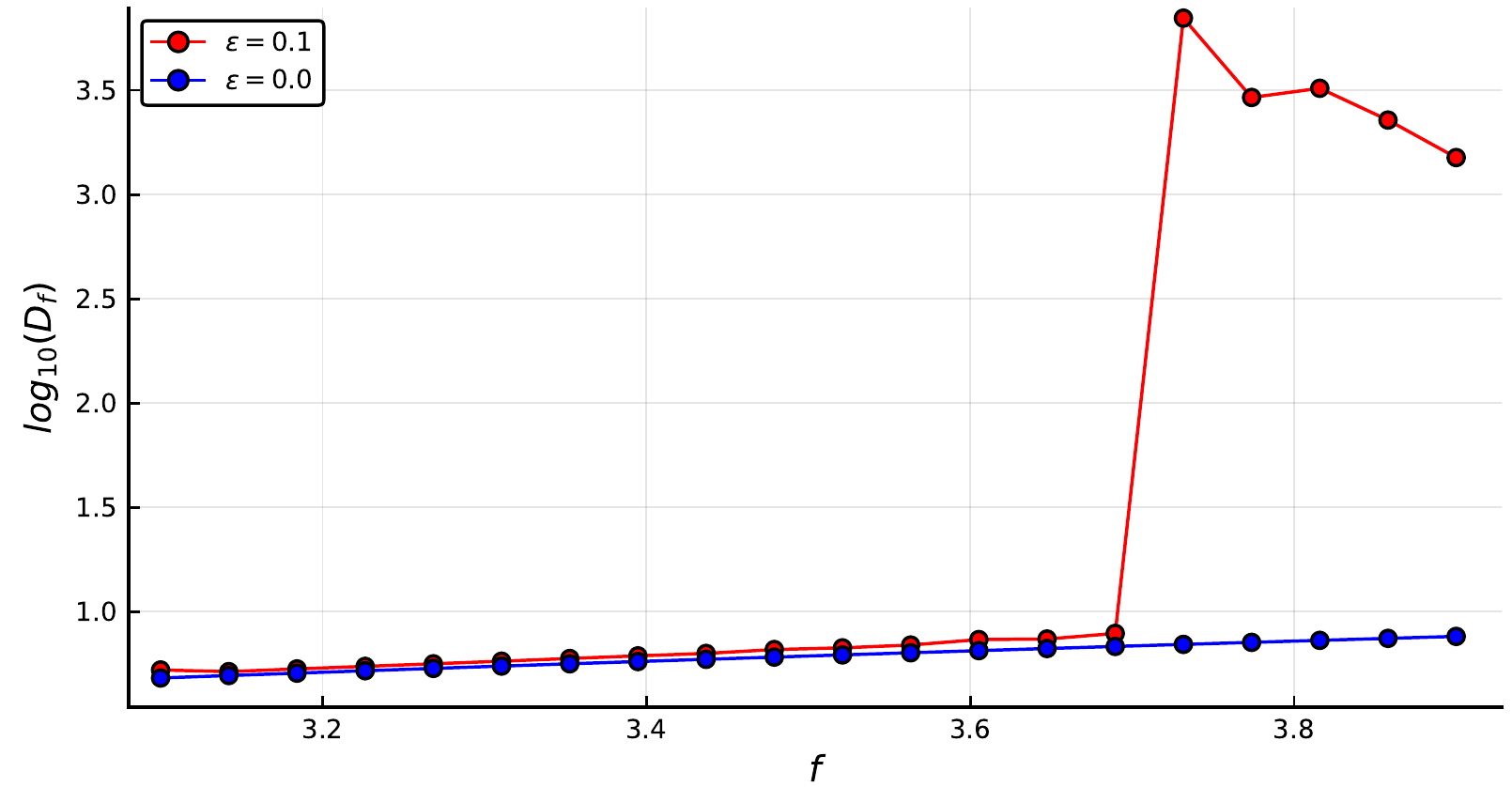} }					
	\caption{Values of the numerical index $\log_{10}\left(\mathcal{D}_f\right)$ [see Eq.~\eqref{average}] as a function of the forcing amplitude $f$ for the original Reid model (blue) and the nonlinear model I (red) systems, showing that the nonlinear model does experience supratransmission at $f_{cr}\approx 3.68$, while the quasilinear one does not. The parameters are set to $\left(N, c, k,\omega\right) = \left(100, 0.1,0.1,3\right)$.\label{fig55}}
\end{figure}

Finally, let us use the proposed numerical criterion to investigate the occurrence of supratransmission for the quasi-linear and nonlinear systems. As we can see in Fig.~\ref{fig55}, for parameter values  $\left(N, c, k,\omega\right) = \left(100, 0.1,0.1,3\right)$, the quasi-linear model proposed by Reid, \textit{fails to exhibit} this phenomenon, since there is no jump at the graph of the numerical criterion index $\log_{10}\left(\mathcal{D}_f\right)$, while the nonlinear model I does show supratransmission. 

\section{Supratransmission and wave packet spreading}

\subsection{Model II: The case of nearest neighbor hysteretic damping}

Note that the nonlinear extension of the periodically forced Reid's model I of Section 2 can be written for $j = 1,\ldots,N$ as follows [see also Eq.~\eqref{eqm}]: 
\begin{equation}\label{eqm1}
M\ddot{x}_j + c\,x_j\, \text{sgn}\left(x_j\dot{x}_j\right) + k\, \left(\Delta^{-}x_j + \Delta^{+}x_j\right) + \epsilon \, \left(\left(\Delta^{-}x_j\right)^3 +\left(\Delta^{+}x_j\right)^3 \right) = 0,
\end{equation}
where we have introduced for brevity of notation the difference operators  $\Delta^{-}x_j := x_j-x_{j-1}$ and $\Delta^{+}x_j := x_j-x_{j+1}$, and similarly, $\Delta^{-}\dot{x}_j := \dot{x}_j-\dot{x}_{j-1}$ and $\Delta^{+}\dot{x}_j := \dot{x}_j-\dot{x}_{j+1}$ for the velocity terms. 

It is now instructive to introduce another model (called model II here) whose hysteretic damping terms are {\it nearest neighbor dependent} and lead to the equations of motion 
\begin{align}\label{eqm2}
M\ddot{x}_j &+ c\,\left(\Delta^{-}x_j\,\text{sgn}(\Delta^{-}x_j \Delta^{-}\dot{x}_j) + c\,\Delta^{+}x_j\,\text{sgn}(\Delta^{+}x_j \Delta^{+}\dot{x}_j)\right) \\ \nonumber 
&+ k\, \left(\Delta^{-}x_j + \Delta^{+}x_j\right) + \epsilon \, \left(\left(\Delta^{-}x_j\right)^3 +\left(\Delta^{+}x_j\right)^3 \right) = 0
\end{align}
Model II is subject to the same boundary conditions as model I [see Eq.~(\ref{eq:bc})], i.e.:
\begin{equation}\label{eqm2dr}
x_0(t) = f\sin{\omega\,t} \,\, \text{ and } \,\, x_{N+1}(t)=0,\,\,\,\, t\in \mathcal{T}\subseteq\mathbb{R}^+
\end{equation}
for the forced case, and
\begin{equation}
x_0(t) = x_{N+1}(t)=0,\,\, t\in \mathcal{T}\subseteq\mathbb{R}^+,
\end{equation}
for the unforced case.

\subsection{Supratransmission in model II}

\begin{figure}[h!]
	\centering
			\resizebox{0.75\columnwidth}{!}{
	\includegraphics{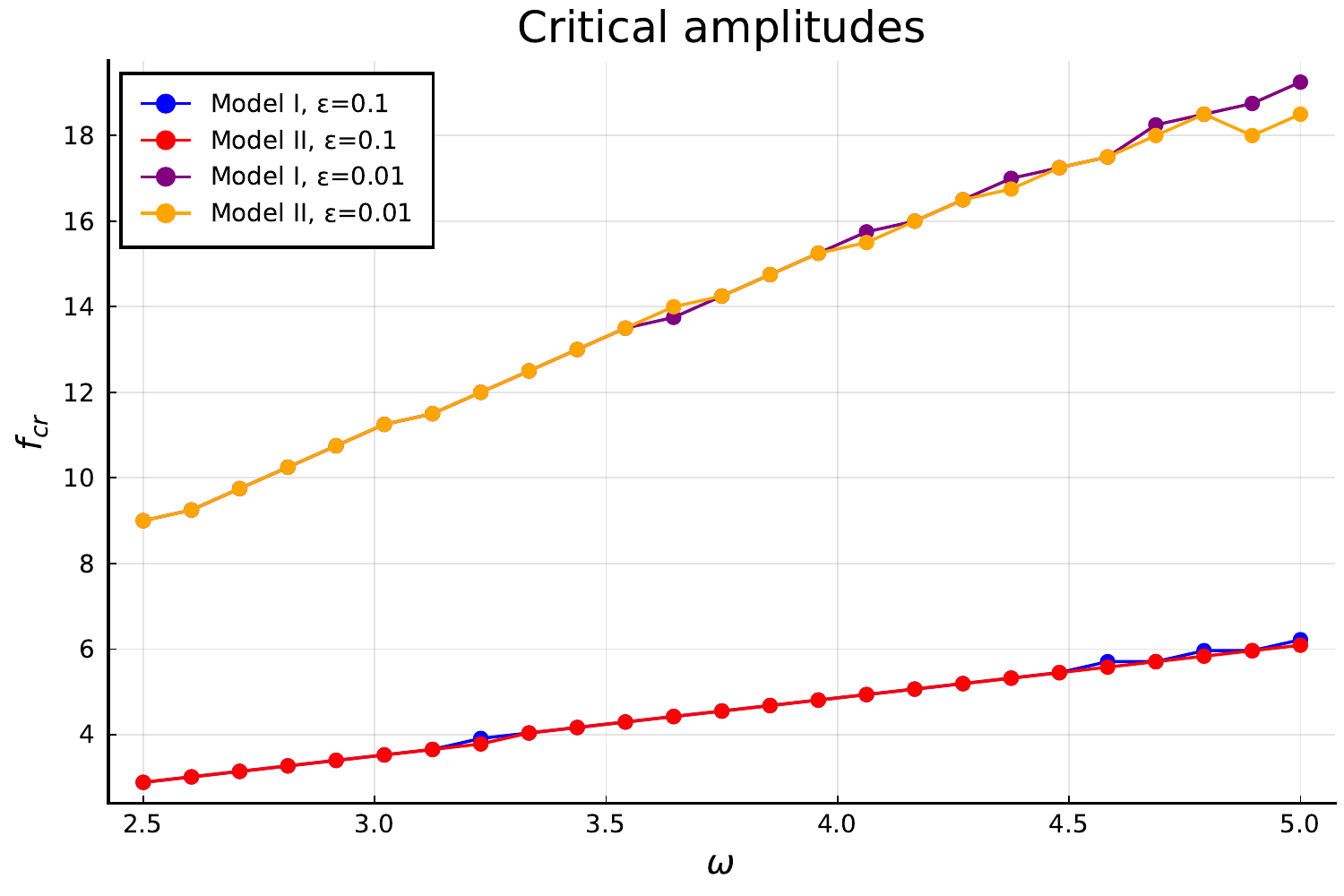} }					
	\caption{Plots of supratransmission critical amplitudes $f_{cr}$ as functions of the forcing frequency $\omega$ for models I and II, for $\epsilon \in \left\{0.01,0.1\right\}$ and $\left(k,c\right) = (0.3, 0.01)$. The curves corresponding to the same $\epsilon$ value practically overlap. Compare with Fig.~\ref{fig9}.}\label{fig4c}
\end{figure}

In this subsection we examine possible differences between our models I and II with regard to their supratransmission properties, plotting in Fig.~\ref{fig4c} for both models the corresponding  critical amplitudes as functions of the forcing frequency $\omega$, with $\epsilon \in \left\{0.01,0.1\right\}$. In these simulations, we have used arrays of $N=200$ coupled oscillators with control parameters $\left(k,c\right) = (0.3, 0.01)$.

It appears that critical amplitudes are nearly identical for the two models, with small discrepancies evident for relative large $\omega$ and $\epsilon$ values.
This implies that the type of hysteretic damping used in the two models does not play a significant role in the emergence of supratransmission phenomena.

\subsection{Wave packet spreading in models I and II}

Another way to compare the dynamics of models I and II is by studying their wave packet spreading properties under time evolution \cite{skokos2018,skokos2020}, for both the forced $(f\neq0)$ and unforced $(f=0)$ cases. To this end, we form an initial wave packet starting with both systems at rest, except for a region of length $L$ in the middle of the lattice where we impose nonzero initial displacements. We then compute the time evolution of the two systems and compare how the wave packets spread in the two models respectively. 

\begin{figure}[h!]
	\centering
			\resizebox{0.9\columnwidth}{!}{
	\includegraphics{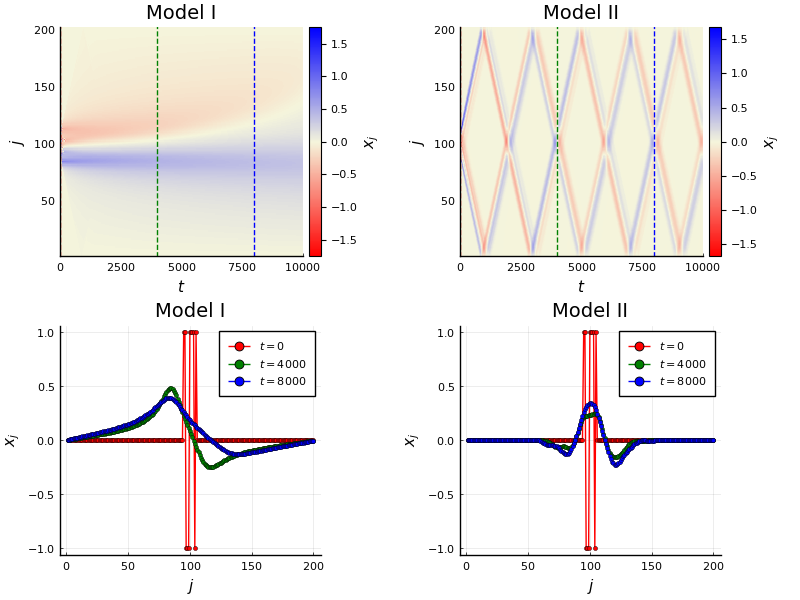} }					
	\caption{Top row: Energy maps of the displacement solutions of the {\it unforced} model I (left) and model II (right), for the parameters $\left(c,k,\epsilon\right)=(0.01, 0.01, 0.1)$. Bottom row: Snapshots of the oscillator displacements $x_j$ versus the index $j$, for time values 0 (red points), 4000 (green points) and 8000 (blue points). The associated times are indicated with vertical lines in the top row energy maps.}\label{fig4d}
\end{figure}

\begin{figure}[h!]
	\centering
	\resizebox{0.9\columnwidth}{!}{
		\includegraphics{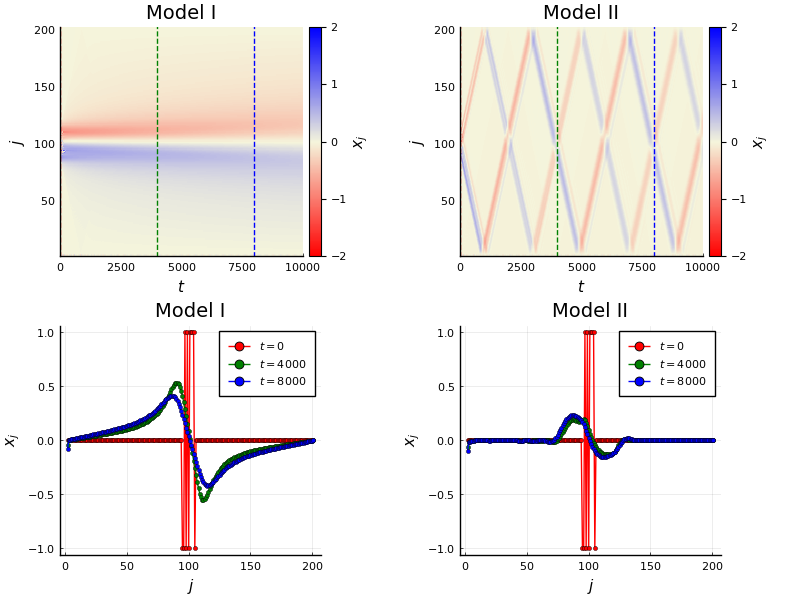} }					
	\caption{Top row: Energy maps of the displacement solutions for the {\it periodically forced} model I (left) and model II (right). The parameters are $\left(c,k,f,\omega,\epsilon\right) = (0.01, 0.01, 2, 3, 0.1)$. Bottom row: Snapshots of the oscillator displacements $x_j$ versus the index $j$, for time values 0 (red points), 4000 (green points) and 8000 (blue points). The associated times are indicated with vertical lines in the top row heatmaps.}\label{fig4f}
\end{figure}

In Figs.~\ref{fig4d} and \ref{fig4f} we present the results obtained from these simulations using arrays of $N=200$ coupled oscillators and total integration time $t=10,000$. Fig.~\ref{fig4d} shows the two wave forms for the {\it unforced}  models initialized with the same initial conditions, i.e.~all particles at rest except for the central packet of $L=11$ sites that are excited randomly at +1 or -1 with equal probabilities. The parameters used are $\left(c,k,\epsilon\right)=(0.01, 0.01, 0.1)$. 

Remarkably, we observe strikingly different behaviors in the the two models: In the unforced case we witness in the top row of Fig.~\ref{fig4d} that model I exhibits a very slow and uniform spreading of the packet in both directions, while model II shows high velocity localized waves that are periodically reflected at the boundaries. We also deduce from different snapshots of the oscillator displacements in the bottom row of Fig.~\ref{fig4d} that model I exhibits a wider range of high displacement sites compared with model II as time increases. \newline

Continuing our comparison of the two models, we plot in Fig.~\ref{fig4f} results obtained in the periodically forced case, using nonzero damping and the same initial excitations for the central $L=11$ sites as before, for the parameters $\left(c,k,f,\omega,\epsilon\right) = (0.01, 0.01, 2, 3, 0.1)$. While we do observe a pattern similar to the corresponding unforced systems of Fig.~\ref{fig4d}, the spreading of the initial excitation in model I appears more uniform, while model II displays out of phase waves propagating in opposite directions and reversing their phase upon collision with the boundaries! Furthermore, we no longer observe a ``shift'' of the maximal displacement sites in time, as in the unforced case of Fig.~\ref{fig4d}.

\section{Stochastic Perturbations of Nonlinear Reid oscillator arrays}

Finally, we proceed to analyze the behavior of arrays of nonlinear Reid oscillators under the effect of random periodic forcing \cite{ciszak,si2007}, by applying additive noise only to the first particle of model I that is periodically driven. Thus, we now solve the same equations of motion as in Eq.~\eqref{eqm}: 
\begin{align}\label{seqm}
M\ddot{x}_j = - c \abs{x_j} \tanh ( \tau \dot{x}_j ) - k\, \left(-x_{j-1} + 2x_j - x_{j+1}\right) - \epsilon \, \left(- \left(x_{j+1}-x_j\right)^3 +\left(x_j-x_{j-1}\right)^3 \right),
\end{align}
for $j=1,\ldots,N$, subject to boundary conditions of the form:
\begin{equation}
x_0(t) = f\sin{\omega\,t} + \sigma\, \xi(t) \,\,\,\,\, \text{ and } \,\,\,\,\, x_{N+1}(t)=0,\,\,\,\, t\in \mathcal{T}\subseteq\mathbb{R}+,
\end{equation}
where $\sigma\, \xi(t)$ represents additive white Gaussian noise and $\sigma$ denotes its intensity (see e.g. \cite{yamgoue2007noise,bodo2009noise}).

Note that this type of random driving transforms the system of (nonlinear) ODEs (\ref{seqm}) to a system of (nonlinear) Stochastic Differential Equations (SDEs). To solve these SDEs, we employ a modified Euler--Heun method with adaptive time step and perform the integration using the DifferentialEquations.jl library \cite{rackauckas2017differentialequations} of the Julia programming language \cite{bezanson2017julia}. 

\begin{figure}[h!]
	\centering
			\resizebox{\columnwidth}{!}{
	\includegraphics{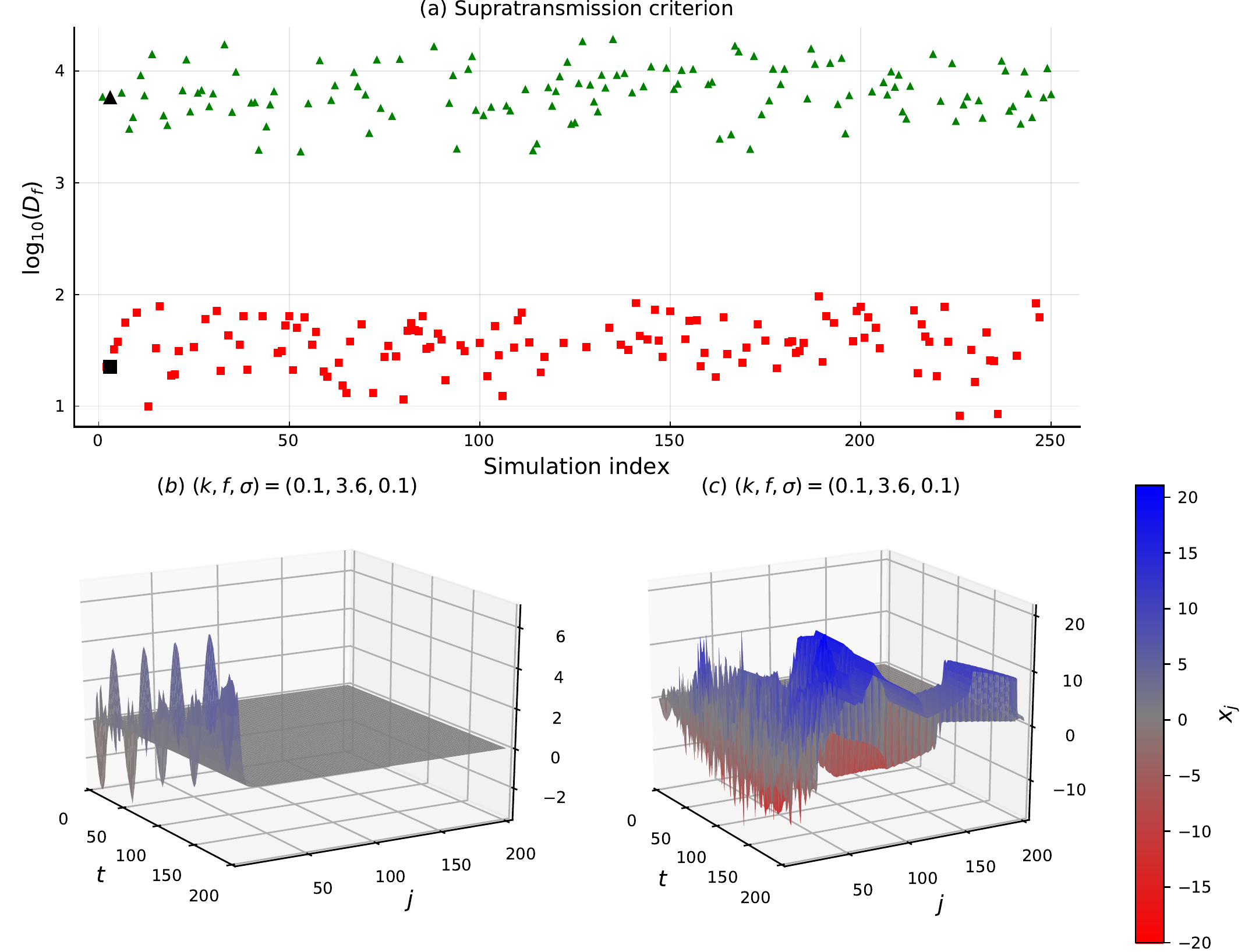} }						
	\caption{$(\alpha)$ Semi-logarithmic plots of 250 simulations of the stochastic system. We present the values of the numerical criterion index $\log_{10}\left(\mathcal{D}_f\right)$, for arrays of size $N = 200$. $(b)-(c)$ Plots of the particle displacements solutions of Eqs.~(\ref{seqm}) for two different realizations of the stochastic system  indicated in $(a)$ by a black square $(b)$ and a black triangle $(c)$. The other parameters are set to $\left(c,\omega,\epsilon\right) = \left(0.1,3,0.1\right)$. \label{fig51}}
\end{figure}

It is particularly interesting that, for a fixed noise intensity, we also observe here the emergence of supratransmission phenomena. This is illustrated in Fig.~\ref{fig51}, where we use our numerical criterion of Eqs.~(\ref{average}), (\ref{critical}), to identify noise realizations that result in supratransmission. More specifically, we set our parameter values at  $\left(c,k,f,\omega,\epsilon,\sigma\right) = \left(0.1, 0.1, 3.6, 3, 0.1, 0.1\right)$ and perform stochastic integration of a system of $N=200$ oscillators up to time $t=200$, for $250$ noise realizations. For each simulation we display in Fig.~\ref{fig51}(a) the value of the numerical criterion $\log_{10}\left(\mathcal{D}_f\right)$ and observe the formation of two clusters. 

The ``lower'' cluster of red squares corresponds to realizations that do not exhibit supratransmission, while the ``upper'' cluster of green triangles corresponds to realizations that do exhibit supratransmission. In Figs.~\ref{fig51}(b),(c) we present indicative plots of the displacement solutions for simulations corresponding to a red square and a green triangle respectively, depicted in black color in Fig.~\ref{fig51}(a). Remarkably, the only difference between the two simulations lies in the realization of the noise process and indicates that for noisy systems there is a probability of supratransmission for values of the forcing amplitude \textit{lower} than the critical amplitude of the corresponding deterministic system! Note that the probability of supratransmission (or lack thereof) is estimated as the ratio of the number of green (or red) triangles over the total number of realizations. 

\begin{figure}[h!]
    \centering
                \resizebox{0.85\columnwidth}{!}{
        \includegraphics{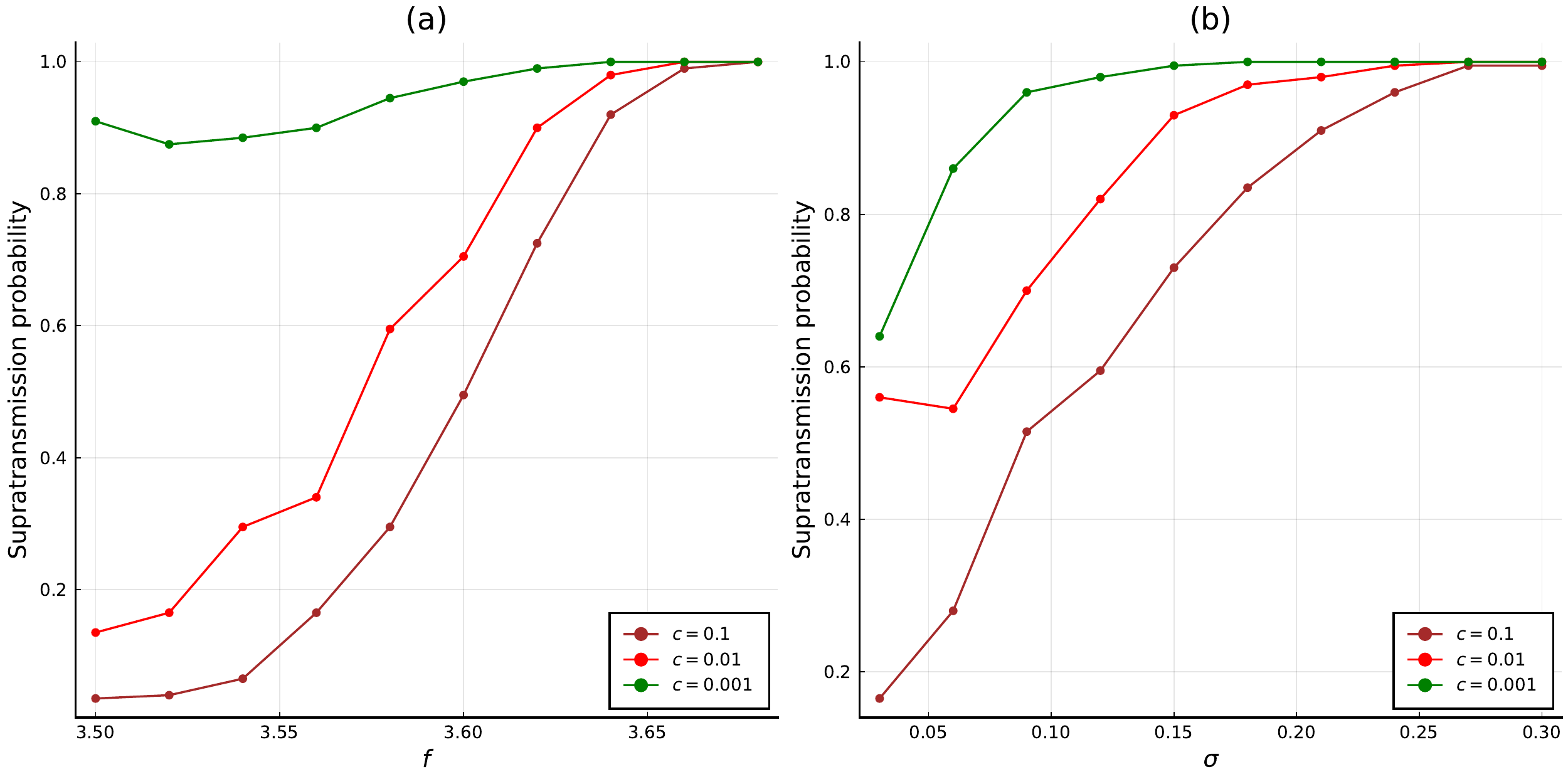} }                                       
    \caption{Probabilities of supratransmission emergence with respect to the forcing amplitude (left) and the noise intensity (right) for varying damping $c \in \left\{0.001, 0.01, 0.1\right\}$. The other parameters are set to $\left(k,\omega,\epsilon\right) = \left(0.1,3,0.1\right)$. \label{fig52}}
\end{figure}

In Fig.~\ref{fig52}, we present additional results regarding the effects of forcing amplitude and noise intensity on the supratransmission probability for parameter values $\left(k,\omega,\epsilon\right) = \left(0.1,3,0.1\right)$, in an array of $N=200$ oscillators. Evidently, higher values of either $f$ or $\sigma$ give rise to higher probabilities of supratransmission phenomena.

\begin{figure}[h!]
	\centering
				\resizebox{\columnwidth}{!}{
	\includegraphics{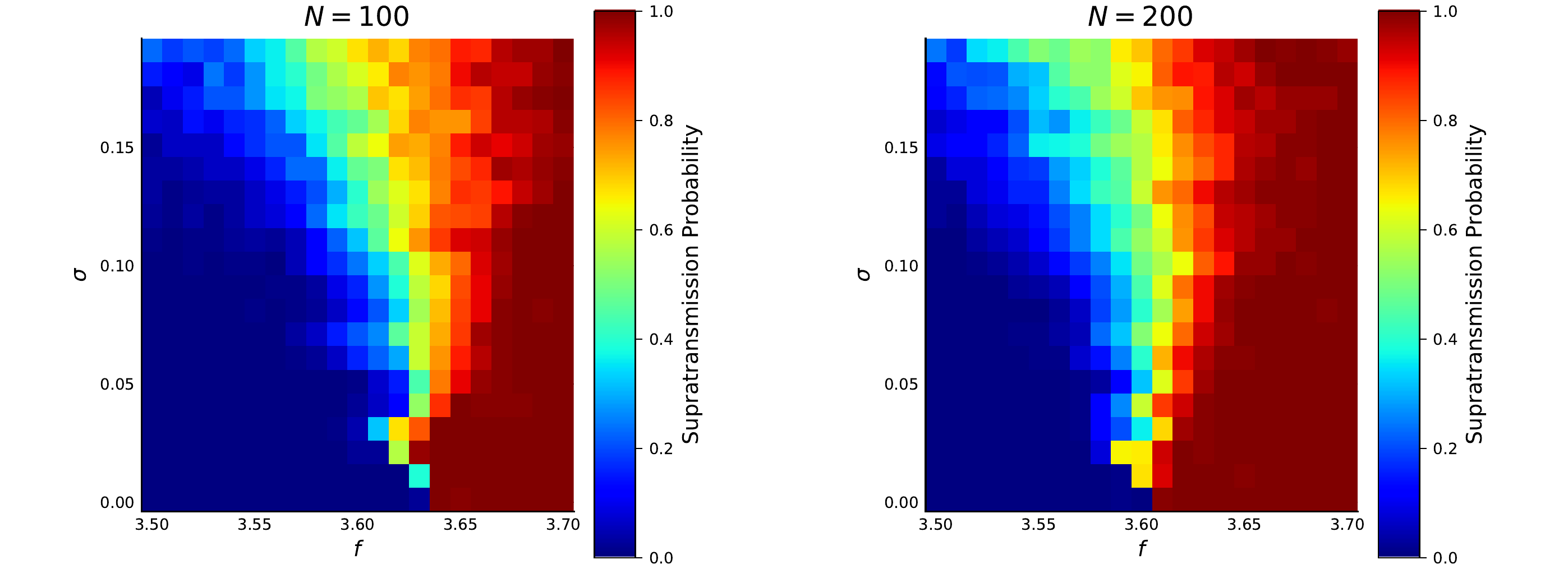} }						
	\caption{Color maps of the supratransmission emergence probability vs. the forcing amplitude $f$ and the noise intensity $\sigma$, for varying array sizes $N \in \left\{100, 200\right\}$. The other parameters are set to $\left(c, k,\omega,\epsilon\right) = \left(0.1,0.1,3,0.1\right)$. 
\label{fig53}}
\end{figure}

Finally, in Fig.~\ref{fig53} we present color maps of the supratransmission probability, taking into account \textit{both} $f$ and $\sigma$, for array sizes $N \in \left\{100, 200\right\}$. At very low values of noise intensity, the system behaves more ``deterministically'' as anticipated, meaning that it has a ``0--1'' like behavior in terms of forcing amplitudes that give rise to supratransmission. As the noise intensity grows, however, the probability of lower amplitudes to generate supratransmission is significantly increased. Thus, we may conclude that the perturbed (noisy) system can excite the required nonlinear modes more ``easily'', i.e.~at lower forcing amplitudes with respect to the associated unperturbed (deterministic) system.

\begin{figure}[h!]
	\centering
					\resizebox{0.85\columnwidth}{!}{
		\includegraphics{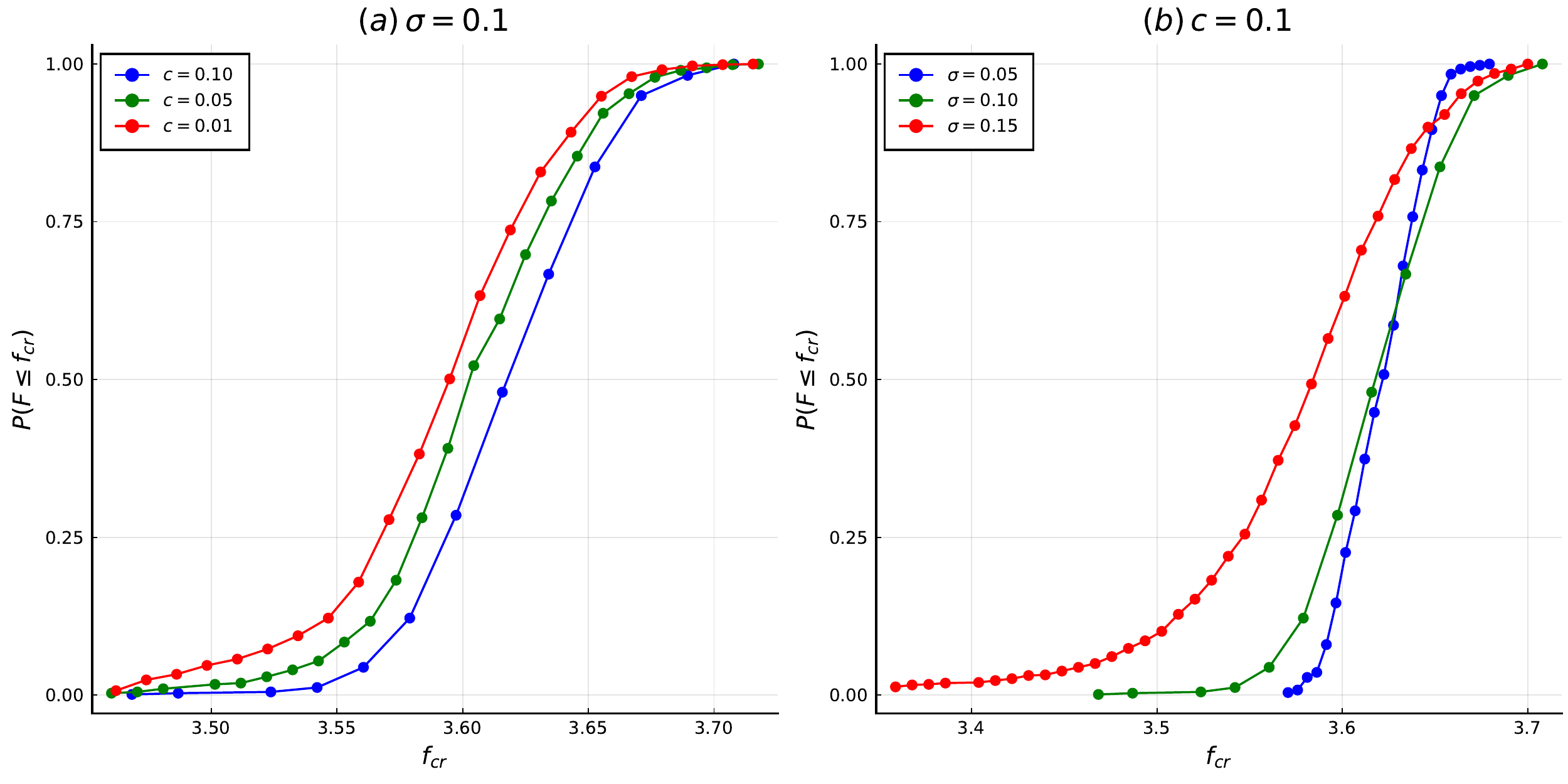} }					
	\caption{Empirical distribution functions of supratransmission critical amplitudes, for varying values of the damping $c$ and noise intensity $\sigma$. The parameters used were $\left(c, k,\omega,\epsilon\right) = \left(0.1,0.1,3,0.1\right)$. \label{fig54}}
\end{figure}

It is also interesting to evaluate the distribution of critical amplitudes for different parameters of the system. To this end, we gradually increase the forcing amplitude until supratransmission occurs, while keeping the same noise realization. In this way, we obtain cumulative distribution functions of $f_{cr}$, as shown in Fig.~\ref{fig54}. Note that the variation of the damping coefficient $c$ in Fig.~\ref{fig54}(a) does not significantly affect the $f_{cr}$ distribution. On the other hand, increasing the noise intensity $\sigma$ leads, as anticipated, to a widening of the support of the distribution and thus to a higher probability that relatively low values of the forcing amplitude may result in supratransmission phenomena.

\section{Traveling breathers and ``breather arrest'' in Models I and II}

Discrete breathers are time-periodic, localized oscillations emerging in spatially extended discrete systems of Hamiltonian lattices (see \cite{FlachGorb} for a comprehensive review). Their occurrence is a result of localized nonlinearities and the discreteness of the lattice. Interestingly, a novel phenomenon termed ``breather arrest'', appearing in lattices of oscillators with {\it viscous damping}, has been analyzed numerically and experimentally in some recent works \cite{Breather3,Breather4}. Essentially, breather arrest refers to the ``birth'' and ``death'' of oscillating breather excitations, following an impulsive excitation at one end of the lattice and an ensuing penetration of these excitations through a small number of oscillators, whose oscillations are damped out eventually after some time.

Here, we apply to both our Models I and II, an impulsive excitation at the left free boundary, with all particles being initially at rest. Specifically, we replace the driving $x_0(t)=fsin(\omega t)$ in our equations (\ref{eq:bc}) and (\ref{eqm2dr}) by $x_0(t)=f\hat{g}(t)$, where $\hat{g}(t)$ is non-zero only over the first half cycle of a sine function with frequency $\omega$. Our right hand boundary is kept fixed at $x_{N+1}=0$. In all our computations we have used $N = 200$ oscillators, parameters $(c, k, f, \omega, \epsilon) = (0.01, 0.3, 3.5, 4, 0.1)$ and integrated up to time $t = 500$, for both Models I and II.
\begin{figure}[h!]
	\centering
				\resizebox{\columnwidth}{!}{
	\includegraphics{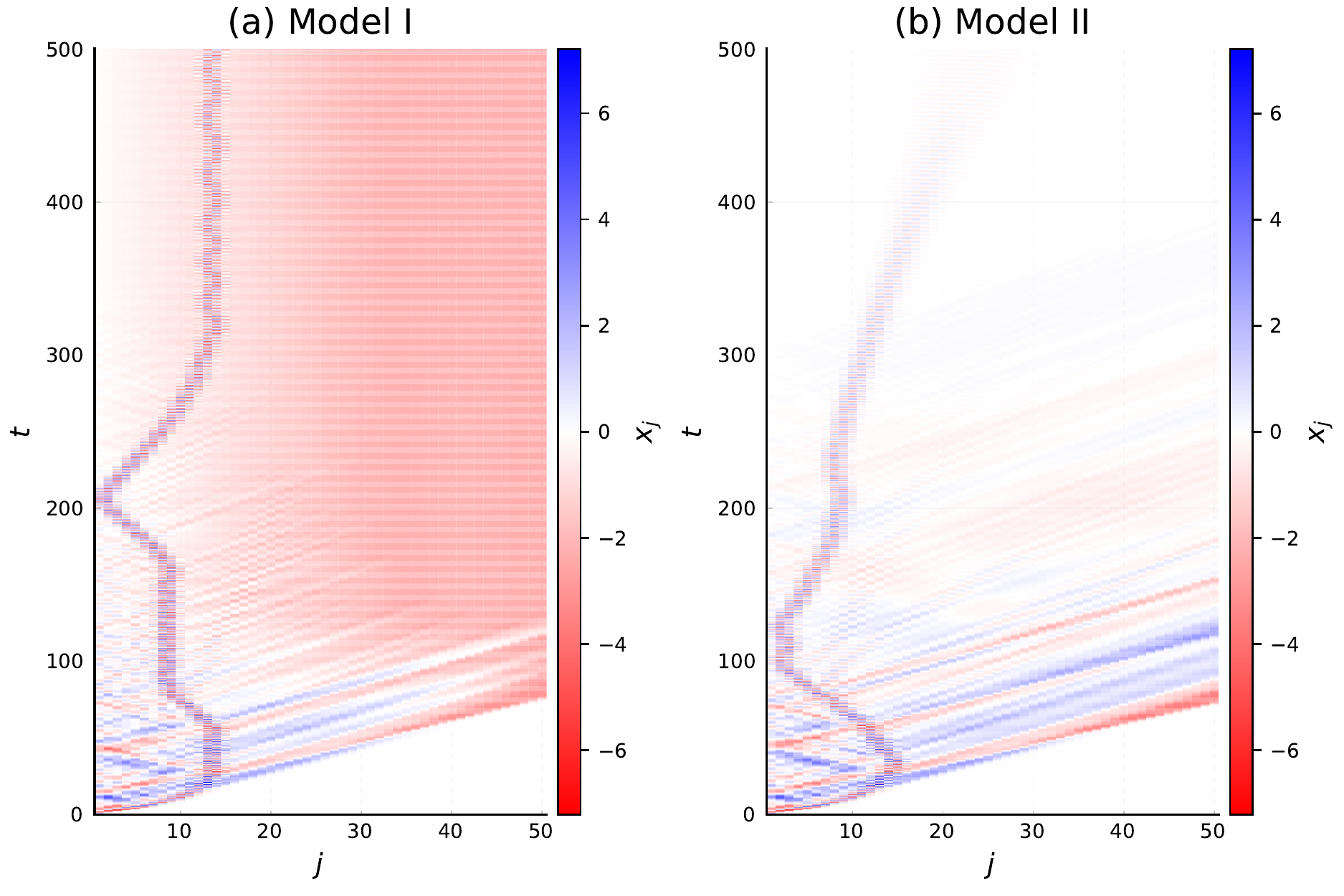} }						
	\caption{Spatiotemporal evolution of the displacement solutions of the impulsively forced Model I (a) and Model II (b). The parameters are set to $(c, k, f, \omega,\epsilon) = (0.01, 0.3, 3.5, 4, 0.1).$ Only the first $50$ of the $N = 200$ particles are depicted. 
\label{fig56}}
\end{figure}
\begin{figure}[h!]
	\centering
				\resizebox{\columnwidth}{!}{
	\includegraphics{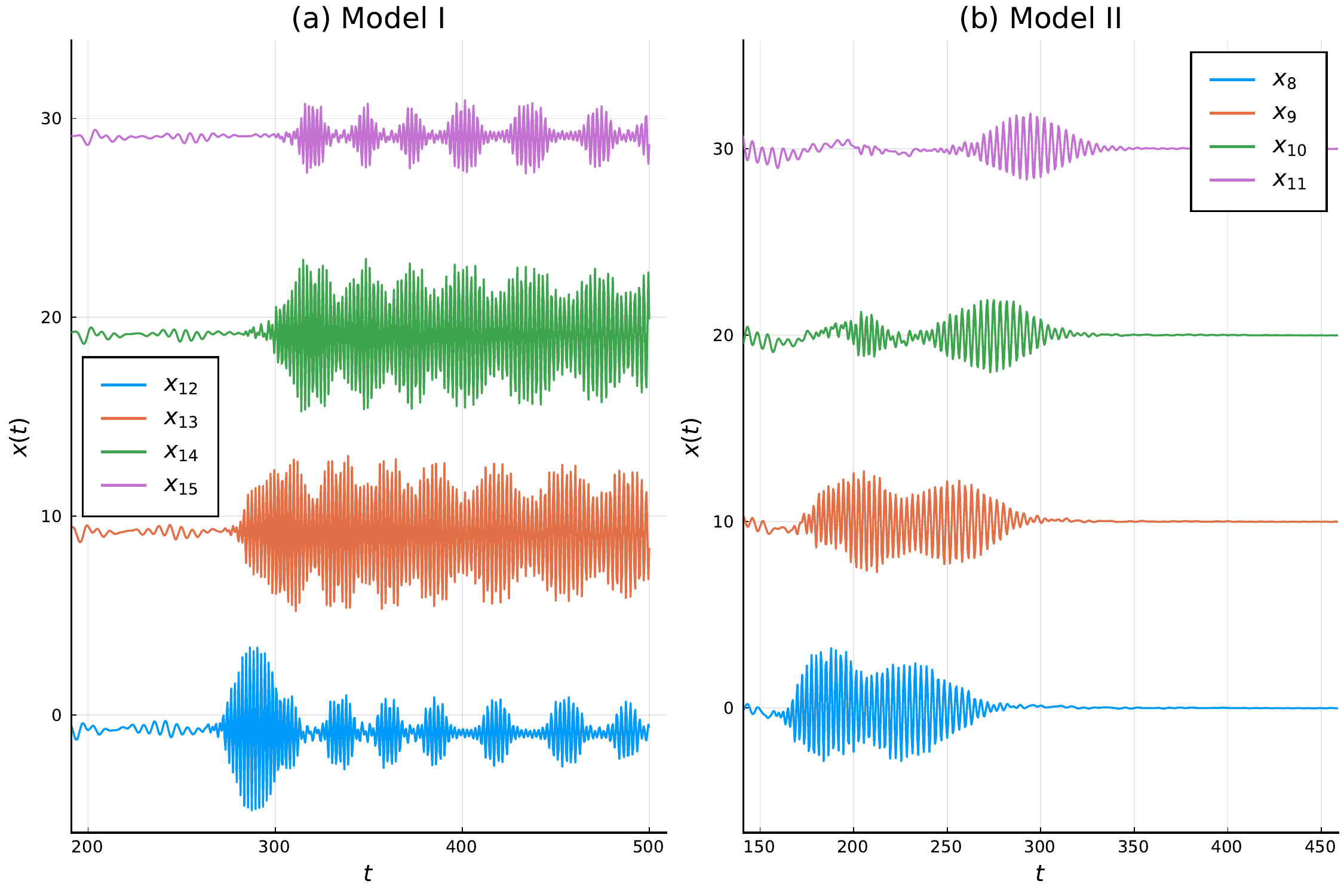} }						
	\caption{Temporal evolution of the displacement solutions of the impulsively forced Model I (a) and model II (b). The parameters are the same as in Fig. \ref{fig56}. We depict only oscillators 12-15 for Model I and 8-11 for Model II, where the breather solutions are clearly visible. We have shifted positions vertically for clarity of presentation. 
\label{fig57}}
\end{figure}

In Fig. \ref{fig56} we present the spatiotemporal evolution of the dynamics of the impulsively forced Models I (a) and II (b). In both cases, we observe the emergence of breathers, but with some notable differences: In Model I, the breather is initially localized near particle 8, but then turns towards the origin and after a reflection at time t = 200, remains localized around oscillator 12. In Model II, the traveling breather is reflected earlier, around t = 100, and then propagates to the right, with rapidly decreasing amplitude.

In Fig. \ref{fig57}, we restrict our attention to a few particles of the two models, and depict the solutions at vertically shifted positions for clarity of presentation. The generated breathers are clearly discernible, having the form of oscillating wavepackets. Clearly, even though both models exhibit oscillations of decreasing amplitude both in time (as well as along the lattice), the phenomenon of ``breather arrest'' in the sense of oscillations dying out after only one maximum \cite{Breather3,Breather4} occurs only in Model II, while Model I exhibits repeated maxima in the oscillations, more like what was observed in \cite{Breather1,Breather2}. Interestingly, in the latter two references, nearest-neighbor coupling was used {\it with viscous damping}, while our Model II employs the same non-local coupling but {\it with hysteretic damping}.

\section{Summary and concluding remarks}

In this paper, we have studied energy transport phenomena associated with supratransmission and wave packet spreading in two different 1-dimensional oscillator arrays with nearest neighbor linear and nonlinear (cubic) stiffness terms in the presence of hysteretic (i.e.~frequency independent) damping. Our main aim was to investigate possible differences in two models, one proposed by Reid with local hysteretic damping \cite{reid1956free} (model I) and one proposed by us (model II), with (non - local) nearest neighbor dependent hysteretic damping.

Our results demonstrate that the system of {\it nonlinearly} coupled Reid oscillators does exhibit the supratransmission property, in contrast to the quasi-linear one, while supertransmission also occurs in model II. In fact, keeping every other parameter fixed leads, in both models to a nearly linear increase of the critical amplitude $f_{cr}$ as a function of the driving frequency $\omega$. In addition, our numerical experiments demonstrate that in both models the damping parameter does not significantly affect the behavior of the associated critical amplitudes. On the other hand, the coupling strength and nonlinearity parameter do play an important role, as {\it higher} values of $k$ or $\epsilon$ lead to {\it lower} values of $f_{cr}$. 

We must stress, of course, that in this paper we have simply studied the onset of supratransmission in our models by only considering driving frequencies $\Omega > 2$, lying outside the harmonic band of the corresponding linear lattices. From the analysis of wave phenomena in similar systems, however, it may turn out that we also need to ensure that our frequencies do lie within the so-called ``nonlinear passband'' of the lattice, otherwise, no propagation will occur. However, the determination of this ``nonlinear passband'' requires quite an intensive analytical and/or numerical investigation, which we prefer to postpone for a future investigation.

We also compared models I and II with regard to their wave packet spreading properties, starting with an initial wave packet located at the center of the array and computing its time evolution within each model. In both forced and unforced cases, we found significantly different behaviors: In particular, in model I a slow uniform spreading away from the center of the lattice is observed in both directions, while model II exhibits rapidly moving localized wave packets that are repeatedly reflected at the boundaries as time evolves.  

We also show that the phenomenon of supratransmission is robust under the effect of additive Gaussian white noise applied to the periodic driving at the first particle of the chain. In fact, the presence of noise perturbations lead to significantly {\it lower} critical amplitudes for supratransmission than in the unforced case. Moreover, as the damping parameter is decreased, a larger noise intensity leads to the widening of the support of the $f_{cr}$ distribution and thus to an even higher probability that lower values of the forcing amplitude can lead to the occurrence of supratransmission.

We also analyzed in our section 5 wave phenomena that occur if we change our periodic to impulsive forcing, as done e.g., in a number of recent papers \cite{Breather1,Breather2,Breather3,Breather4}, where the authors observed the phenomenon of ``breather'' formation and its degradation, termed ``breather'' arrest, under viscous damping. We were thus able to compare our results with those found in the literature and found many similarities but also some differences, showing e.g., that our Model II with nearest-neighbor coupling exhibits the more commonly observed phenomenon of single breather propagation rather than decaying quasiperiodic oscillations of multi--breather structure.

Concerning engineering applications, it might be possible to investigate energy transport in material ``strings'' of nonlinear oscillators and use our results to identify the particular type of damping and/or nonlinear potential terms that best fit the response of the ``string'' to periodic forcing at one of its ends. In this regard, it would be interesting to replace the ODEs employed here by {\it fractional} differential equations, so as to account for non-locality in the range of particle interactions
(see e.g. \cite{bountis2019energy,macias2021nonlinear}.

Regarding a theoretical explanation of our results, it would be interesting to derive a PDE in the continuum limit of our models and show that it possesses nonlinear wave solutions that become unstable at the onset of supratransmission. Thus, it might be possible to replace our numerical criterion with one related to the variation of the total energy of the system (see e.g. \cite{dauxois2007modulational}) and investigate the possible noise--induced emergence of breather modes within parameter ranges that do not support such modes in the noise--free case \cite{bodo2009noise}.

\subsection*{Acknowledgments}
We are indebted to the referees for their very useful remarks, which helped us considerably improve our paper. TB acknowledges that the results of Sections 2 and 5 were obtained under the scientific project No. 21-71-30011 of the Russian Science Foundation, and  those of 3.1 and 3.2 under the Grant No. AP08856381 of the Science Committee of the Ministry of Education and Science of the Republic of Kazakhstan, KazNU, Institute of Mathematics and Mathematical Modeling. CSp acknowledges partial support for this work by funds from the Ministry of Education and Science of Kazakhstan, in the context of the Nazarbayev University internal grant ``Rapid response fixed astronomical telescope for gamma ray burst observation (RARE)'' (OPCRP2020002).
%
%

\end{document}